\documentclass{article}
\begin{document}
\centerline{\textbf{Ramanujan type $1/\pi$ Approximation Formulas}}
\[
\]
\centerline{\bf Nikos Bagis}\\
\centerline{Stenimahou 5 Edessa Pellas 58200, Greece.}
\centerline{e-mail: nikosbagis@hotmail.gr}
\[
\]
\begin{quote}
\begin{abstract}
In this article we use theoretical and numerical methods to evaluate in a closed-exact form the parameters of Ramanujan type $1/\pi$ formulas. 
\end{abstract}

\bf keywords \rm{$\pi$-formulas; Ramanujan; elliptic functions; singular modulus; alternative modular bases; approximations; numerical methods}
\end{quote}

\section{Introduction}

We give the definitions of the Elliptic Integrals of the first and second kind respectively (see [9],[4]):
\begin{equation}
K(x)=\int^{\pi/2}_{0}\frac{dt}{\sqrt{1-x^2\sin^2(t)}} \textrm{ and } E(x)=\int^{\pi/2}_{0}\sqrt{1-x^2\sin^2(t)}dt . 
\end{equation}
In the notation of Mathematica we have
\begin{equation}
K(x)=\textrm{EllipticK}[x^2] \textrm{ and } E(x)=\textrm{EllipticE}[x^2] . 
\end{equation}
Also we have (see [9],[7]):  
\begin{equation}
\dot{K}(k)=\frac{dK(k)}{dk}=\frac{E(k)}{k(1-k^2)}-\frac{K(k)}{k} . 
\end{equation}
The elliptic singular moduli is defined to be the solution of the equation: 
\begin{equation}
\frac{K\left(\sqrt{1-w^2}\right)}{K(w)}=\sqrt{r} .
\end{equation}
In Mathematica is stated as
\begin{equation}
w=k=k_r=k[r]=\textrm{InverseEllipticNomeQ}[e^{-\pi \sqrt{r}}]^{1/2} . 
\end{equation}
The complementary modulus is given by $k^{'2}_{r}=1-k_r^2$.\\  
Also we will need the following relation of the elliptic alpha function (see [7]):
\begin{equation}
a(r)=\frac{\pi}{4K(k_r)^2}-\sqrt{r}\left(\frac{E(k_r)}{K(k_r)}-1\right) . 
\end{equation}
The Hypergeometric functions are defined by
$$
{}_{m+1}F_{m}\left(a_1,a_2,\ldots,a_{m+1};b_1,b_2,\ldots,b_m;z\right):=$$
\begin{equation}
=\sum^{\infty}_{n=0}\frac{(a_1)_n(a_2)_n\ldots(a_{m+1})_n}{(b_1)_n(b_2)_n\ldots(b_m)_n}\frac{z^n}{n!} , \textrm{ for } |z|<1 , 
\end{equation}
and $(a)_0:=1$, $(a)_n:=a(a+1)(a+2)\ldots(a+n-1)$, for each positive integer $n$.

\section{The construction of some $1/\pi$ and $1/\pi^2$ formulas}

It holds 
\begin{equation}
\phi_1(z)={}_3F_{2}\left(\frac{1}{2},\frac{1}{2},\frac{1}{2};{1,1};z\right)=\frac{4K^2\left(\frac{1}{2}(1-\sqrt{1-z})\right)}{\pi^2} . 
\end{equation} 
Consider the following equation with respect to the function $\phi_1(z)$:
$$
\sum^{\infty}_{n=1}\frac{\left(\frac{1}{2}\right)^3_n}{(n!)^3}z^n(an+b)=\frac{g}{\pi}\Leftrightarrow b\phi_1(z)+az\phi'_1(z)=\frac{g}{\pi} .   
$$
Set $w=1/2\left(1-\sqrt{1-k^2}\right)$, $1-2w=\sqrt{1-z}=k'_r$.\\
But
$$
b\phi_1(z)+az\phi'_1(z)=\frac{g}{\pi}\Leftrightarrow g=\frac{4K(w)(aE(w)+(b+a(w-1)-2bw)k(w))}{\pi(1-2w)} , 
$$
hence
$$
\sum^{\infty}_{n=1}\frac{\left(\frac{1}{2}\right)^3_n}{(n!)^3}4^n(w-w^2)^n(an+b)=$$
$$
=\frac{4K(\sqrt{w})\left(a E(\sqrt{w})+(b-a+aw-2bw-2bw)K(\sqrt{w})\right)}{\pi^2(1-2w)} . 
$$
For $w=k_r$ we get
$$
\sum^{\infty}_{n=1}\frac{\left(\frac{1}{2}\right)^3_n}{(n!)^3}4^n(k_rk'_r)^{2n}(an+b)=$$
\begin{equation}
=\frac{4K(k_r)\left(a E(k_r)+(b-a+ak^2_r-2bw-2bk^2_{r})K(k_r)\right)}{\pi^2(1-2k^2_r)} . 
\end{equation}
Now using the formula for $a(r)$, in the sense that
\begin{equation}
E(k_r)=K(k_r)-\frac{a(r)K(k_r)}{\sqrt{r}}+\frac{\pi}{4K(k_r)\sqrt{r}} , 
\end{equation}
for suitable values for $a$, $b$, $c$ we get the following theorem:
\[
\]
\textbf{Theorem 2.1}
\begin{equation}
\sum^{\infty}_{n=1}\frac{\left(\frac{1}{2}\right)^3_n}{(n!)^3}4^n(k_rk'_r)^{2n}\left(\sqrt{r}(1-2k^2_r)n+a(r)-\sqrt{r}k^2_r\right)=\frac{1}{\pi}
\end{equation}
\[
\]
\textbf{Example.}
$$
\sum^{\infty}_{n=0}\frac{\left(\frac{1}{2}\right)^3_n}{(n!)^3}(40\sqrt{2}-56)^n(an+b)=\frac{4a}{7\pi}+\frac{5a}{7\sqrt{2}\pi}+4(-4a+\sqrt{2}a+14b)\frac{\Gamma^2\left(\frac{9}{8}\right)}{7\pi\Gamma^2\left(\frac{5}{8}\right)} . 
$$
From which a special case is
$$
\sum^{\infty}_{n=0}\frac{\left(\frac{1}{2}\right)^3_n}{(n!)^3}(40\sqrt{2}-56)^n(n+\frac{2}{7}-\frac{1}{7\sqrt{2}})=\frac{8+5\sqrt{2}}{14\pi} . 
$$
\[
\]
\textbf{Theorem 2.2}
\begin{equation}
\sum^{\infty}_{n=0}\frac{B^{(2)}_n}{(n!)^2}(k_r)^{2n}(\sqrt{r}k'^2_r n+a(r)-\sqrt{r}k_r^2)
=\frac{1}{\pi} . 
\end{equation}
\textbf{Proof.}\\
We use the function
\begin{equation}
\phi_2(z)={}_2F_1\left(\frac{1}{2},\frac{1}{2};1;z\right)=\frac{2K(\sqrt{z})}{\pi} . 
\end{equation}
Then
if
$$ B^{(2)}_n:=\sum^{n}_{j=0}\left[\left(^n_j\right)\left(\frac{1}{2}\right)_n\left(\frac{1}{2}\right)_{n-j}\right]^2
$$
\begin{equation}
\phi^2_2(z)=...=\sum^{\infty}_{n=1}\frac{z^n}{(n!)^2}\sum^{n}_{j=0}\left[\left(^n_j\right)\left(\frac{1}{2}\right)_n\left(\frac{1}{2}\right)_{n-j}\right]^2 , 
\end{equation}
where
$$
c\phi_2(z)+bz\phi'_2(z)+az^2\phi''_2(z)=
\sum^{\infty}_{n=0}\frac{B^{(2)}_n}{(n!)^2}z^{n}(an^2+(b-a)n+c)
$$
Hence we get
$$
\sum^{\infty}_{n=0}\frac{B^{(2)}_n}{(n!)^2}k_r^{2n}(an^2+(b-a)n+c)
=\frac{2\left(a E(k_r)+(2b-2bk_r^2-4a+6ak_r^2)E(k_r) K(k_r)\right)}{\pi^2(1-k_r^2)^2}+
$$
$$
+\frac{2(3a-2b+2c+(-4a+2b-2c)k_r^2)K(k_r)}{\pi^2 (1-k_r^2)} . 
$$
For $a=0$, $b=1$, $c=(-k_r^2+a(r)r^{-1/2})k'^{-2}_r$, we get 
\[
\]
\textbf{Theorem 2.3} Set
\begin{equation}
B^{(3)}_n:=\sum^{n}_{j=0}\left[\left(^n_j\right)\left(\frac{1}{2}\right)_n\left(\frac{1}{2}\right)_{n-j}\right]^3 , 
\end{equation}
then an $1/\pi^2$ formula is the following
\begin{equation}
\sum^{\infty}_{n=0}\frac{B^{(3)}_n}{(n!)^3}(2k_rk'_r)^{2n}(n^2+(b(r)-1)n+c(r))=\frac{3}{(1-2k_r^2)^2 r \pi^2}
\end{equation}
where
$$
b(r)=\frac{3a_r+\sqrt{r}-6a(r)k_r^2-9\sqrt{r}k_r^2+12\sqrt{r}k_r^4}{\sqrt{r}(1-2k_r^2)^2}
$$
and
$$
c(r)=\frac{3a(r)^2-6a(r)\sqrt{r}k_r^2-rk_r^2+4rk_r^4}{r(1-2k_r^2)^2}
$$
\textbf{Proof.}\\
Set 
$$
\phi_3(z)={}_3F_2\left(\frac{1}{2},\frac{1}{2},\frac{1}{2};1,1;z\right)^2=\left(\frac{16K^2\left(\frac{1}{2}(1-\sqrt{1-z})\right)}{\pi^2}\right)^2 , 
$$
then 
$$
c\phi_3(z)+bz\phi'_3(z)+az^2\phi''_3(z)=
\sum^{\infty}_{n=0}\frac{B^{(3)}_n}{(n!)^3}z^{n}(an^2+(b-a)n+c)
$$
The left hand of the above equation is a function of $E(x)$, $K(x)$, and can evaluated when we set certain values to the parameters $a$, $b$, $c$.
\[
\]
\textbf{Examples.}\\ 
\textbf{1)}
$$
\frac{1}{1200(161\sqrt{5}-360)\pi^2}=
$$
\begin{equation}
\sum^{\infty}_{n=0}\frac{B^{(3)}_n}{(n!)^3}\left(51841-23184\sqrt{5}\right)^{n}\left(n^2+\left(1-\frac{521}{288\sqrt{5}}\right)n+\frac{5}{12}-\frac{521}{576\sqrt{5}}\right)
\end{equation}
\textbf{2)}
\footnotesize

$$
b(163)=\frac{191211325848427}{151931373056001}-\frac{1010784962625383717350772720\cdot 2^{2/3}}{151931373056001 \left(B_1-\sqrt{489} B_2\right)^{1/3}}-\frac{4\cdot 2^{1/3} \left((B_1-\sqrt{489}B_2\right)^{1/3}}{151931373056001}
$$

$$
B_1=5680848001702137216093843898647314524189
$$

$$
B_2=76896989960589381643149203281167
$$

$$-5839006481108705728+9529627071955041072\cdot b(163)-4530513053635162884\cdot b(163)^2+$$
$$+668649972819460401\cdot b(163)^3=0$$

$$
c(163)=\frac{14178679829869760}{24764813808128163}-\frac{4 \left(C_1-\sqrt{489} C_2\right)^{1/3}}{24764813808128163}-\frac{6241484569597616793758909818952\cdot 2^{2/3}}{24764813808128163 \left(C_3-\sqrt{489} C_4\right)^{1/3}}
$$

$$
C_1=5512985602111283751597893407219881834715037026
$$

$$
C_2=101526256966667546381077303112958296550
$$

$$
C_3=2756492801055641875798946703609940917357518513
$$

$$
C_4=50763128483333773190538651556479148275
$$

$$
-24380823840878077184+13131020889593608594752\cdot c(163)-$$
$$-30513780896384581928640\cdot c(163)^2+17765361127840243394169\cdot c(163)^3=0
$$

\normalsize

\begin{equation}
\sum^{\infty}_{n=0}\frac{4^nB^{(3)}_n}{(n!)^3}(k_{163}k'_{163})^{2n}(n^2+(b(163)-1)n+c(163))=\frac{A}{\pi^2}
\end{equation}
\footnotesize

$$
A=\frac{4 \left(12660947754667+26680 \left(A_1-\sqrt{489} A_2\right)^{1/3}+26680 \left(A_1+\sqrt{489} A_2\right)^{1/3}\right)}{8254937936042721}
$$

$$
A_1=106866398697613339845357037
$$

$$
A_2=3136555671686449089
$$

$$
y_{163}=(k_{163}k'_{163})^2=\frac{1}{16}-\frac{266933400}{\left(-1+557403 \sqrt{489}\right)^{1/3}}+\frac{10005}{2} \left(-1+557403 \sqrt{489}\right)^{1/3}
$$

$$
-1+16408588290048048\cdot y_{163}-768\cdot y_{163}^2+4096\cdot y_{163}^3=0
$$
\normalsize
\[
\]
Formula (18) gives about 17 digits per term and is a formula for $1/\pi^2$. For $r=253$ we have another such formula which gives 21 digits per term constructed in the same way as (18).

\section{The study of a non usual $1/\pi$ formula}

The $j$ invariant is given by (see [17]):
\begin{equation}
j(z)=\left(\left(\frac{\eta(z/2)}{\eta(z)}\right)^{16}+16\left(\frac{\eta(z)}{\eta(z/2)}\right)^8\right)^3 , 
\end{equation}
where $z=\sqrt{-r}$, $r$-positive real and
$$
\eta(z)=e^{\pi i z/12}\prod^{\infty}_{n=1}\left(1-e^{2\pi i n z}\right)
$$
is the Dedekind eta function.\\ 
Also 
\begin{equation}
\frac{\eta(z)}{\eta(z/2)}=\frac{k^{1/12}_r}{2^{1/6}k'^{1/6}_r} . 
\end{equation}
From [24] section 7, Theorem 7.4 and from [11] formula (5.8), when $q=e^{2\pi i z}$, $z=\sqrt{-r}$, $r$ positive real, the modular $j$-invariant is also given by
\begin{equation}
j(z)=1728\frac{Q^3(q)}{Q^3(q)-R^2(q)}.  
\end{equation}
where 
$$
P(q)=1-24\sum^{\infty}_{n=1}\frac{nq^n}{1-q^n} \textrm{ , } Q(q)=1+240\sum^{\infty}_{n=1}\frac{n^3q^n}{1-q^n}$$
and 
$$
R(q)=1-504\sum^{\infty}_{n=1}\frac{n^5q^n}{1-q^n} .
$$
The function $t_r$ is given from
\begin{equation}
t_r=\frac{Q_r}{R_r}\left(P_r-\frac{6}{\pi\sqrt{r}}\right) ,
\end{equation}
where
$$P_r=P(-e^{-\pi \sqrt{r}})\textrm{ , } Q_r=Q(-e^{-\pi\sqrt{r}})\textrm{ and } R_r=R(-e^{-\pi\sqrt{r}}) . $$

i) Using Theorems 3 and 4 of [25], relation (21) equivalently can be transformed to 
\begin{equation}
j(z)=\frac{432}{\beta_{r}(1-\beta_{r})} .
\end{equation}
Also note that we have
\begin{equation}
j(\sqrt{-r})=j_r=\frac{256(1-k^2_r+k^4_r)^3}{(k_rk'_r)^4}=\frac{432}{\beta_r(1-\beta_r)} . 
\end{equation}
Hence with our method in [25] we can simplify the known results of [24] and [11] using the function $\beta_r$, which defined as the root of the equation: 
\begin{equation}
\frac{{}_{2}F_{1}\left(\frac{1}{6},\frac{5}{6};1;1-w\right)}{{}_{2}F_{1}\left(\frac{1}{6},\frac{5}{6};1;w\right)}=\sqrt{r} . 
\end{equation}  

ii) Set now $m_r:=k_r^2$ and let $a(r)$, $E(x)$ be the elliptic alpha function and the complete elliptic integral of the second kind respectively (see [7],[4]), then:  
$$t_r=\frac{1}{(1-2\beta_{r/4})u_{r/4}^2}\left(P(q)-\frac{6}{\sqrt{r}\pi}\right)=$$
$$
=\frac{1}{(1-2\beta_{r/4})u_{r/4}^2}\left(3\frac{E(m_{r/4})}{K(m_{r/4})}-2+m_{r/4}-\frac{3\pi}{4\sqrt{r/4}K(m_{r/4})^2}\right)F_{r/4}^2
$$
or
\begin{equation}
t_r=\frac{1+m_{r/4}-\frac{6}{\sqrt{r}}a\left(\frac{r}{4}\right)}{\sqrt{1-m_{r/4}+m_{r/4}^2}(1-2\beta_{r/4})} .
\end{equation}
Hence from the above evaluations and the $1/\pi$ series in [6] and [11] we get the next reformulation:
\[
\]
\textbf{Theorem 3.1}
If we define 
\begin{equation}
J_r:=1728j_r^{-1}=4\beta_{r}(1-\beta_{r})
\end{equation} 
\begin{equation}
T_r:=\frac{1+k_{r}^2-\frac{3}{\sqrt{r}}a\left(r\right)}{\sqrt{1-k^2_r+k^4_r}(1-2\beta_{r})}=\frac{2j^{1/3}_r\sigma\left(r\right)G_r^{8}}{\sqrt{r}\sqrt{j_r-1728}}
\end{equation} 
then
\begin{equation}
\frac{3}{\pi\sqrt{r}\sqrt{1-J_{r}}}=\sum^{\infty}_{n=0}\frac{\left(\frac{1}{6}\right)_n\left(\frac{5}{6}\right)_n\left(\frac{1}{2}\right)_n}{(n!)^3}(J_{r})^n(6n+1-T_r)
\end{equation}
\[
\]
\textbf{Note.}
The function $G_r$ is the Weber invariant and  $$\sigma(r)=2\sqrt{r}(1+k^2_r)-6a(r)$$ (see [7],[5] chapter 5).\\
The above formulas (27), (28) and (29) can be used for numerical and theoretical evaluations.     
\[
\]
\centerline{\textbf{Similarities of formula (29) and a fifth order base formula}}
\[
\]
From the identity 
\begin{equation}
{}_2F_{1}\left(\frac{1}{6},\frac{5}{6};1;\frac{1-\sqrt{1-z}}{2}\right)^2={}_3F_{2}\left(\frac{1}{6},\frac{5}{6},\frac{1}{2};1,1;z\right) , 
\end{equation}
and using the following relations found in [7]:  
\begin{equation}
K_s(x)=\frac{\pi}{2}{}_2F_1\left(\frac{1}{2}-s,\frac{1}{2}+s;1;x^2\right) \textrm{ and } E_s(x)=\frac{\pi}{2}{}_2F_1\left(-\frac{1}{2}-s,\frac{1}{2}+s;1;x^2\right)
\end{equation}
\begin{equation}
E_s=(1-k^2)K_s+\frac{k(1-k^2)}{1+2s}\dot{K}_s \textrm{ , } \dot{K}_s(t)=\frac{dK_s(t)}{dt}
\end{equation}
\begin{equation}
a_s(x_r):=\frac{\pi}{4K_s(x_r)}\frac{\cos(\pi s)}{1+2s}-\sqrt{r}\left(\frac{E_s(x_r)}{K_s(x_r)}-1\right) , 
\end{equation}
 with $s=1/3$ one can get, (working as in Theorem 2.1) the following Ramanujan-type $1/\pi$ formula: 
\begin{equation}
\sum^{\infty}_{n=0}\frac{\left(\frac{1}{6}\right)_n\left(\frac{5}{6}\right)_n\left(\frac{1}{2}\right)_n}{(n!)^3} (4\beta_r(1-\beta_r))^n)\left(3n-5\frac{\beta_r-\frac{a_5(r)}{\sqrt{r}}}{1-2\beta_r}\right)=\frac{3}{2\pi\sqrt{r}(1-2\beta_r)} ,  
\end{equation}
where the function $\alpha_5(r)=a_{1/3}(\sqrt{\beta_r})$ is algebraic for $r\in \bf Q^{*}_{+}$.\\ 
The parameters and the corresponding function $\alpha_5(r)$ of (34) are those of fifth singular moduli base theory. Also (34) in comparison with (29) gives the following theorem. 
\[
\]
\textbf{Theorem 3.2}
\begin{equation}
10\alpha_5(r)r^{-1/2}=10a_{1/3}(\sqrt{\beta_r})r^{-1/2}=1+8\beta_r-\frac{1+k^2_r-3a(r)r^{-1/2}}{\sqrt{1-k^2_r+k^4_r}}
\end{equation}  

The above formula is for general evaluation of elliptic alpha function in the fifth elliptic base.
\[
\]
Also from the cubic theory as in fifth, we have 

\begin{equation}
{}_3F_{2}\left(\frac{1}{3},\frac{2}{3},\frac{1}{2};1,1;w\right)={}_2F_{1}\left(\frac{1}{3},\frac{2}{3};1;\frac{1-\sqrt{1-w}}{2}\right)^2
\end{equation}
we get
\begin{equation}
\sum^{\infty}_{n=0}\frac{\left(\frac{1}{3}\right)_n\left(\frac{2}{3}\right)_n\left(\frac{1}{2}\right)_n}{(n!)^3}[4\alpha_3(r)-4\alpha^2_3(r)]^n(n-b)=\frac{\sqrt{3}}{2\pi\sqrt{r}(1-2\alpha_3(r))}
\end{equation}

\begin{equation}
b=\frac{4\left(\alpha_3(r)-a_{1/6}[\alpha_3^{1/3}(r)]r^{-1/2}\right)}{3(1-2\alpha_3(r))}
\end{equation}

\section{Examples and Evaluations}

1)For $r=2$
$$
J_2=\frac{27}{125}
$$
$$
T_2=\frac{5}{14}
$$
and
\begin{equation}
\frac{15\sqrt{5}}{14\pi}=\sum^{\infty}_{n=0}\frac{\left(\frac{1}{6}\right)_n\left(\frac{5}{6}\right)_n\left(\frac{1}{2}\right)_n}{(n!)^3}\left(\frac{27}{125}\right)^n\left(6n+\frac{9}{14}\right)
\end{equation}
2) For $r=4$ we have
\begin{equation}
\frac{11 \sqrt{\frac{11}{3}}}{14 \pi }=\sum^{\infty}_{n=0}\frac{\left(\frac{1}{6}\right)_n\left(\frac{5}{6}\right)_n\left(\frac{1}{2}\right)_n}{(n!)^3}\left(\frac{8}{1331}\right)^n\left(6n+\frac{10}{21}\right)
\end{equation}
3) For $r=5$ we have
$$T_5=\frac{1}{418} \left(139+45 \sqrt{5}\right)$$
$$J_5=\frac{27 \left(-1975+884 \sqrt{5}\right)}{33275}$$
Hence
$$
\frac{\sqrt{21650+5967 \sqrt{5}}}{\pi }=
$$
\begin{equation}
=\sum^{\infty}_{n=0}\frac{\left(\frac{1}{6}\right)_n\left(\frac{5}{6}\right)_n\left(\frac{1}{2}\right)_n}{(n!)^3}\left(\frac{-53325+23868 \sqrt{5}}{33275}\right)^n\left(836 n+93-15 \sqrt{5}\right) 
\end{equation}
4) For $r=8$ we have 
$$k_8^2=113+80 \sqrt{2}-4 \sqrt{2 \left(799+565 \sqrt{2}\right)}$$
$$a(8)=2 \left(10+7 \sqrt{2}\right) \left(1-\sqrt{-2+2 \sqrt{2}}\right)^2
$$
Then
$$
\frac{15 \sqrt{\frac{5}{2} \left(84125+81432 \sqrt{2}\right)}}{9982 \pi }=
$$
\begin{equation}
=\sum^{\infty}_{n=0}\frac{\left(\frac{1}{6}\right)_n\left(\frac{5}{6}\right)_n\left(\frac{1}{2}\right)_n}{(n!)^3}\left(\frac{5643000-3990168 \sqrt{2}}{1520875}\right)^n\left(\frac{3276-1125 \sqrt{2}+29946 n}{4991}\right)
\end{equation}
5) For $r=18$ we have 
$$k_{18}=(-7+5\sqrt{2})(7-4\sqrt{3})$$
$$a(18)=-3057+2163 \sqrt{2}+1764 \sqrt{3}-1248 \sqrt{6}$$
$$\alpha_{6}=\frac{1}{500}(68-27\sqrt{6})$$
$$\beta_{18}=\frac{1}{2}-\frac{7 \left(49982+4077 \sqrt{6}\right)}{10 \sqrt{5} \left(989+54 \sqrt{6}\right)^{3/2}}$$
\begin{equation}
J_{18}=\frac{637326171-260186472 \sqrt{6}}{453870144125}
\end{equation}
\begin{equation}
T_{18}=\frac{712075+49230 \sqrt{6}}{1074514}
\end{equation}
Hence we get the formula giving  8 digits per term:\\(Note that the number of digits per term is determined by the value of $J_r$, approximately.)
$$ 
\frac{5 \sqrt{23124123365-13274820 \sqrt{6}}}{1074514 \pi }=
$$
\begin{equation}
\sum^{\infty}_{n=0}\frac{\left(\frac{1}{6}\right)_n\left(\frac{5}{6}\right)_n\left(\frac{1}{2}\right)_n}{(n!)^3}\left(\frac{637326171-260186472 \sqrt{6}}{453870144125}\right)^n\times$$
$$\times\left(6n+\frac{9 \left(40271-5470 \sqrt{6}\right)}{1074514}\right)
\end{equation}
6) For $r=27$ $$k_{27}=\frac{1}{2}\sqrt{\frac{1+100\cdot2^{1/3}-80\cdot 2^{2/3}}{2+\sqrt{3-100\cdot2^{1/3}+80\cdot2^{2/3}}}}$$
$$a(27)=3\left[\frac{1}{2}\left(\sqrt{3}+1\right)-2^{1/3}\right]$$
$a(27)$ is obtained from [7] page 172.
$$
J_{27}=\frac{56143116+157058640\cdot2^{1/3}-160025472\cdot 2^{2/3}}{817400375}
$$
$$T_{27}=\frac{58871825+22512960\cdot2^{1/3}+13208820\cdot 2^{2/3}}{132566687}$$
Hence we get the 11 digits per term formula:
$$
\frac{935}{\pi} \sqrt{\frac{935}{3\left(761257259-157058640\sqrt[3]{2}+160025472 \sqrt[3]{4}\right)}}=
$$
$$
=\sum^{\infty}_{n=0}\frac{\left(\frac{1}{6}\right)_n\left(\frac{5}{6}\right)_n\left(\frac{1}{2}\right)_n}{(n!)^3}\left(\frac{56143116+157058640 \sqrt[3]{2}-160025472 \sqrt[3]{4}}{817400375}\right)^n\times
$$
$$
\times\left(6n+\frac{6 \left(12282477-3752160 \sqrt[3]{2}-2201470\sqrt[3]{4}\right)}{132566687}\right)
$$
7) From the Wolfram pages 'Elliptic Lambda Function' and 'Elliptic Singular Value' we have: 
$$
k_{58}=\left(-1+\sqrt{2}\right)^6 \left(-99+13 \sqrt{58}\right)
$$
and
$$
a(58)=\frac{1}{64} \left(-70+99 \sqrt{2}-13 \sqrt{29}\right) \left(5+\sqrt{29}\right)^6 \left(-444+99 \sqrt{29}\right)
$$
Also using the cubic theta identities, (see [25] relations (2),(3),(4),(30)) we evaluate $\alpha_{174}$ numerically to 1500 digits and then $\beta_{58}$ to 1500 digits accuracy. We  then apply the 'Recognize' routine of Mathematica. The result is the minimum polynomial of $\beta_{58}$ (this can be done also from (19) and (23)):  
$$
1-1399837865393267000 x+79684665286353732299517000 x^2-$$
$$-159369327773031733812500000 x^3+79684663886515866906250000 x^4=0 . 
$$
Solving this equation with respect to $x$ we get the value of $\beta_{58}$ in radicals. Thus
\begin{equation}
J_{58}=\frac{1399837865393267-259943365786104 \sqrt{29}}{39842331943257933453125}
\end{equation}
\begin{equation}
T_{58}=\frac{5 \left(1684967251+24160612 \sqrt{29}\right)}{10376469642}
\end{equation}
The result is the formula
$$
\frac{5 \sqrt{\frac{5}{87} \left(13826969809210107-90211316 \sqrt{29}\right)}}{357809298 \pi }=
$$
$$
\sum _{n=0}^{\infty } \frac{\left(\frac{1}{6}\right)_n \left(\frac{5}{6}\right)_n\left(\frac{1}{2}\right)_n }{(n!)^3}\left(\frac{1399837865393267-259943365786104 \sqrt{29}}{39842331943257933453125}\right)^n
$$
\begin{equation}
\times\left(\frac{6117973}{32528118}-\frac{8628790}{25557807 \sqrt{29}}+6 n\right)
\end{equation}
\newline
which gives 18 digits per term.
\[
\]
8) For $r=93$ (see [7] pg.158), we have 
$$
\sigma(93)=6G_{93}^{-6}\left(\frac{\sqrt{3}+1}{2}\right)^3\left(15\sqrt{93}+13\sqrt{31}+201\sqrt{3}+217\right) .
$$
From [5] chapter 34 we have
$$
G_{93}=\frac{\left(3 \sqrt{3}+\sqrt{31}\right)^{1/4} \left(39+7 \sqrt{31}\right)^{1/6}}{2^{1/3}}
$$
also 
$$
a(r)=\sqrt{r}\frac{1+k^2_r}{3}-\frac{\sigma(r)}{6}
$$
$$
G^{-24}_{93}=4k^2_{93}(1-k^2_{93})
$$
Hence
\footnotesize
$$
(k_{93}k'_{93})^2=\frac{1}{224589314596+129666700800 \sqrt{3}+40337431680 \sqrt{31}+23288826960 \sqrt{93}}
$$
$$
J^{-1}_{93}=119562334956358303022500+21474029280866147440000 \sqrt{31}+
$$
$$
+470106000 \sqrt{129368095019778762513344107725+23235195778655878514048710848 \sqrt{31}}
$$
$$
T_{93}=\frac{10559116299575+1317692448000 \sqrt{3}+275805228680 \sqrt{31}-81807235875 \sqrt{93}}{15081520900138} ,
$$

\normalsize
where
$$
j_{93}=1728J^{-1}_{93}
$$
This result is a very flexible formula that gives about 24 digits per term.

\section{Neat Examples with Mathematica and Simplicity}

The class number $h(-d)$, $d\in\bf N\rm$ of the equivalent quadratic forms is given by
\begin{equation}
h(-d)=-\frac{w(d)}{2d}\sum^{d-1}_{n=1}\left(\frac{-d}{n}\right)n , 
\end{equation}
where
$w(3)=6$, $w(4)=4$ else $w(d)=2$.
 $\left(\frac{n}{m}\right)$, is the Jacobi symbol. Observe that $h(-163)=1$ (see [17]). For small values of $h(-d)$ we have greater possibility to evaluate $J_d$ and $T_d$ in radicals.  

The simplest way to evaluate the parameters $J_{163}$ and $T_{163}$ is again with Mathematica.
\[
\]
The general algorithm is:

i) Set $r=d$ and $k[r]=\textrm{InverseEllipticNomeQ}[e^{-\pi\sqrt{r}}]^{1/2}$, then we can evaluate $\beta_r$ and $j_r$ from relations (19) and (23). Hence we get the value of $J_r$ as in section 4 example 7.
   
ii) For the evaluation of $T_r$ we will need the value of $a(r)$ which is given from (see [7]):
\begin{equation}
a(r)=\frac{\pi}{4K^2}-\sqrt{r}\left(\frac{E}{K}-1\right) .
\end{equation}
This in Mathematica is given from
\begin{equation}
a(r)=\frac{\pi}{\textrm{4EllipticK}[k[r]^2]^2}-\sqrt{r}\left(\frac{\textrm{EllipticE}[k[r]^2]}{\textrm{EllipticK}[k[r]^2]}-1\right)
\end{equation}
Hence taking the package
$$<< \textrm{NumberTheory}`\textrm{Recognize}`$$
and
$$\textrm{Recognize}[N[J_{163},1500],16,x]$$
$$\textrm{Recognize}[N[T_{163},1500],16,x]$$
we get two equations. After solving them we get if $r\in \bf N\rm$ (here $r=163$), 
the values of the parameters $J_r$ and $T_r$ in algebraic-closed forms. The results are the $\pi$ formulas.
\[
\]
1) We have that $J_{163}$ is root of
\[
\]
\footnotesize
$$-64+2552810853189232588558727380998000 x-2198253790246041723377943360187500 x^2+$$
$$+224451422498574115473590775022822688001953125 x^3=0$$
\normalsize
hence
\footnotesize
$$
J_{163}=4\frac{C_1-C_2\left(-A_1+\sqrt{489} B_1\right)^{-1/3}+30591288 \left(-A_1+\sqrt{489} B_1\right)^{1/3}}{10792555251621895860488211571345343375} 
$$
\[
\]
\footnotesize
$$
A_1=12737965652562547164590026038483234248161827096523072256574968383
$$
$$
B_1=229038073182066825378006485964950394558349727761749294205546402325349
$$
$$
C_1=8808429913332498766352891
$$
$$
C_2=902206261147132595923169636910570558029813352485594880
$$
\normalsize
\[
\]
From $J_{r}=4 \beta_r (1-\beta_r)$, we get the value of $\beta_r$  
and hence
\footnotesize
$$T_{163}=5\frac{12948195754365757115+8 \left(A_2-B_2 \sqrt{489}\right)^{1/3}+8 \left(A_2+B_2 \sqrt{489}\right)^{1/3}}{83470787671093501833}$$
\normalsize
where
\footnotesize
$$
A_2=3802386862487392962897493239274992371253057854289262
$$ 
$$
B_2=3865464212119923579732688315287754932290919450 
$$
\normalsize
\[
\] 
The above parameters give 32 digits per term 
\[
\]
2) Another evaluation is taking $d=r=253$:
\footnotesize
\[
\]
$$
J_{253}=\frac{A_1-A_2\sqrt{11}+31990140 \sqrt{A_3-A_4\sqrt{11} }}{A_5}
$$
\[
\]
$$A_1=2804365789259959094417576921792857440357087269234369$$
$$A_2=845548099807651569627713349319558464492321957799872$$
$$A_3=1433462642401972199773341051748172965440271797713951$$
$$6818782945906676740858207407330990565$$
$$A_4=43220524871261259540733172862370537466134334936322822$$
$$33926553935879770457716659641968088$$
$$A_5=1066755353338783886372226117351012749877681799897625$$
\normalsize
and
\footnotesize
\[
\]
$$
T_{253}=\frac{1875 \sqrt{B_1-B_2 \sqrt{11}}+3847208393012364625+752271279708923520 \sqrt{11}}{6969874104047710086}
$$
\[
\]
$$B_1=213216899528167866600672118125$$
$$B_2=60533150139616794053500831192$$
\[
\]
\normalsize
The above parameters give 41 digits per term.
\[
\]
\textbf{Conclusion}\\
We have given a way of how we can construct a very large number of Ramanujan's type $1/\pi$ formulas. It is true that in most cases, from $r=1$ to $100$ (or higher), using Mathematica program, such formulas are very simple, as long as $h(-d)$ remains small and the parameters are solutions of solvable polynomial equations.
\[
\]

\centerline{\bf References}\vskip .2in

[1]: M.Abramowitz and I.A.Stegun: Handbook of Mathematical Functions. Dover Publications. (1972).

[2]: B.C.Berndt: Ramanujan`s Notebooks Part I. Springer Verlag, New York. (1985).

[3]: B.C.Berndt: Ramanujan`s Notebooks Part II. Springer Verlag, New York. (1989).

[4]: B.C.Berndt: Ramanujan`s Notebooks Part III. Springer Verlag, New York. (1991).

[5]: B.C. Berndt: Ramanujan's Notebooks Part V. Springer Verlag, New York, Inc. (1998) 

[6]: Bruce C. Berndt and Heng Huat Chan: Ramanujan and the Modular j-Invariant. Canad. Math. Bull. Vol.42(4), (1999). pp.427-440.

[7]: J.M. Borwein and P.B. Borwein: Pi and the AGM. John Wiley and Sons, Inc. New York, Chichester, Brisbane, Toronto, Singapore. (1987). 

[8]: I.S. Gradshteyn and I.M. Ryzhik: Table of Integrals, Series and Products. Academic Press. (1980).
 
[9]: E.T. Whittaker and G.N. Watson: A course on Modern Analysis. Cambridge U.P. (1927)

[10]: I.J. Zucker: The summation of series of hyperbolic functions. SIAM J. Math. Ana.10.192. (1979)

[11]: Bruce.C. Berndt and Heng Huat Chan: Eisenstein Series and Approximations to $\pi$. Page stored in the Web.

[12]: S. Ramanujan: Modular equations and approximations to $\pi$. Quart. J. Math.(Oxford). 45, 350-372. (1914).

[13]: S. Chowla: Series for $1/K$ and $1/K^2$. J. Lond. Math. Soc. 3, 9-12. (1928)

[14]: N.D. Baruah, B.C. Berndt and H.H. Chan: Ramanujan's series for $1/\pi$: A survey. American Mathematical Monthly 116, 567-587. (2009)

[15]: T. Apostol: Modular Functions and Dirichlet Series in Number Theory. Springer

[16]: Bruce.C. Berndt, S. Bhargava and F.G. Garvan: Ramanujan's Theories of Elliptic Functions to Alternative Bases. Transactions of the American Mathematical Society. 347, 4163-4244. (1995)

[17]: D. Broadhurst: Solutions by radicals at Singular Values $k_N$ from New Class Invariants for $N\equiv3\;\; mod\;\; 8$'. arXiv:0807.2976 (math-phy).

[18]: J.V. Armitage W.F. Eberlein: Elliptic Functions. Cambridge University Press. (2006)

[19]: N.D. Baruah, B.C. Berndt: Eisenstein series and Ramanujan-type series for $1/\pi$. Ramanujan J.23. (2010) 17-44

[20] N.D. Baruah, B.C. Berndt: Ramanujan series for $1/\pi$ arising from his cubic and quartic theories of elliptic functions. J. Math. Anal. Appl. 341. (2008) 357-371 

[21]: B.C. Berndt: Ramanujan's theory of Theta-functions. In Theta functions: from the classical to the modern Editor: Maruti Ram Murty, American Mathematical Society. 1993 

[22]: J.M. Borwein and P.B. Borwein: A cubic counterpart of Jacobi's identity and the AGM. Transactions of the American Mathematical Society, 323, No.2, (Feb 1991), 691-701 

[23]: Habib Muzaffar and Kenneth S. Williams: Evaluation of Complete Elliptic Integrals of the first kind at Singular Moduli. Taiwanese Journal of Mathematics, Vol. 10, No. 6, pp 1633-1660, December 2006 

[24]: Bruce C. Berndt and Aa Ja Yee: Ramanujans Contributions to Eisenstein Series, Especially in his Lost Notebook. (page stored in the Web). 

[25]: Nikos Bagis: Eisenstein Series, Alternative Modular Bases and Approximations of $1/\pi$. arXiv:1011.3496 (2010)

\end{document}